\documentclass{amsart}
\title{Differential transcendence and walks on self-similar graphs}
\author{Yakob Kahane}
\address{\'Ecole Polytechnique, Palaiseau, France}
\author{Marni Mishna}
\address{Department of Mathematics, Simon Fraser University, Burnaby, Canada}

\usepackage[utf8]{inputenc}
\usepackage{amsxtra, mathtools}
\usepackage{float}
\usepackage{geometry}
\usepackage{tikz}
\usepackage{enumitem}
\usepackage{nicefrac}
\usepackage{bbm}
\usepackage{xspace}
\usepackage{url}
\usepackage{graphics,color, soul}
\usepackage{graphicx}
\graphicspath{{Images/}}
\usepackage{wrapfig}
\usepackage{subfigure}
\usepackage{fullpage}
\newtheorem{theorem}{Theorem}
\newtheorem*{nontheorem}{Theorem}
\newtheorem{lemma}[theorem]{Lemma}
\newtheorem{corollary}[theorem]{Corollary}
\newtheorem{conjecture}[theorem]{Conjecture}
\newtheorem{prop}[theorem]{Proposition}
\theoremstyle{definition}
\newtheorem{definition}[theorem]{Definition}

\newtheorem{example}{Example}

\newcommand{\C}{\mathbb{C}}
\newcommand{\origin}{\mathbf{o}}
\newcommand{\cell}{\ensuremath{C}}
\newcommand{\hcell}{\ensuremath{\hat\cell}}

\usepackage{mathpazo}

\usepackage{comment}
\usepackage[textwidth=3cm,textsize=tiny]{todonotes}
\newcounter{todocounter}
\newcommand{\mjm}[2] 
{{\color{red}#2}\stepcounter{todocounter}\todo[,color=yellow!25]{\sf\thetodocounter:  #1}}
\newcommand{\yk}[1]
{\stepcounter{todocounter}\todo[color=blue!20]{\tiny \sf\thetodocounter:  #1}}

\begin{document}
\begin{abstract}Symmetrically self-similar graphs are an important type of fractal graph. Their Green's functions satisfy order one iterative functional equations.  We show when the branching number of a generating cell is two, either the graph is a star consisting of finitely many one-sided lines meeting at an origin vertex,  in which case the Green's function is algebraic, or the Green's function is differentially transcendental over $\mathbb{C}(z)$.  The proof strategy relies on a recent work of Di Vizio, Fernandes and Mishna.  The result adds evidence to a conjecture of Krön and Teufl about spectra of the difference Laplacian of this family of graphs.

\end{abstract}
\maketitle


\section{Introduction}
\label{sec:intro}
The discretization of fractals as graphs opens the possibility of using
combinatorial  methods to study self-similarity. In
particular, some types of fractal graphs possess a regularity which can be
exploited in enumerative analysis. The focus here is walks on graphs in the family
of \emph{symmetrically self-similar}\/ fractal graphs. 

A uniform random walk on a graph $X=(V(X), E(X))$ starts at some
vertex, and proceeds along the edges of the graph. We assume the
degree of every vertex is bounded and that every edge incident to a
given vertex is equally likely to be chosen.
We define the following transition matrix $P$ associated to a graph:
For vertices $x,y$ from $V(X)$, the $[x,y]$ entry of~$P$ is the
transition probability
\begin{equation}\label{eq:P}
  p(x,y):=
  \begin{cases}\frac{1}{\deg(x)} &{\{x,y\}\in E(X)}\\0&\text{otherwise}
  \end{cases}.\end{equation}
A random walk on $X$ is a Markov chain
$(V_i)_{i=0}^n$, $V_i\in V(X)$. The probability of the path $V_0=x_0, V_1=x_1, \dots, V_n=x_n$, given that $V_0=x_0$, is
\[
\mathbb{P}_{x_0}[V_0=x, V_1=x_1, \dots, V_n=x_n]=p(x_0, x_1) p(x_1, x_2)\cdots p(x_{n-1}, x_n). 
\]
We denote the probability of a walk that starts at $x$ ends at $y$
after $n$ steps ($\mathbb{P}_x[V_n=y]$ ) by~$p^{(n)}(x,y)$. A \emph{Green's function} of $(X, P)$, denoted $G(x,y|z)$, is a generating function of probabilities: 
\[
G(x, y| z) := \sum_{n\geq 0} p^{(n)}(x,y)\, z^n.
\]
When we mean walks that start and end at a specified origin vertex
labeled $\origin$ we refer to them as \emph{the} Green's function for a graph, and use the shorthand $G(z)$. Precisely,
\[G(z):= G(\origin,\origin|z)=\sum_{n\geq 0} p^{(n)}(\origin, \origin)\, z^n.\] 

This Green's function, and the functional equations it satisfies, can offer significant
information about~$X$, and possibly even objects the
graph might encode (like groups in the case of a Cayley graph).  We focus
here on different notions of transcendency. Recall, a series in the
ring of formal power series $\C[[z]]$ is said to be \emph{algebraic}
if it satisfies a non-trivial polynomial equation, and is
\emph{D-finite} over $\C(x)$ if its derivatives span a finite
dimensional vector space over~$\C(x)$. A function is
\emph{differentiably algebraic} if it satisfies a non-trivial
polynomial differential equation, and is \emph{differentially
  transcendental\footnote{Also called hypertranscendental in some
    literature}} if is not differentiably algebraic.  These categories
have demonstrated themselves to be useful in a variety of
combinatorial contexts. One related example is the case of Cayley graphs: algebraicity and D-finiteness of a Green's function are each correlated with structural properties of the group~\cite{bell_complexity_2020}. 

\begin{figure}
\centering
\noindent\begin{tikzpicture} \node[inner sep=0pt]  at (0,0) {\includegraphics[width=.7\textwidth]{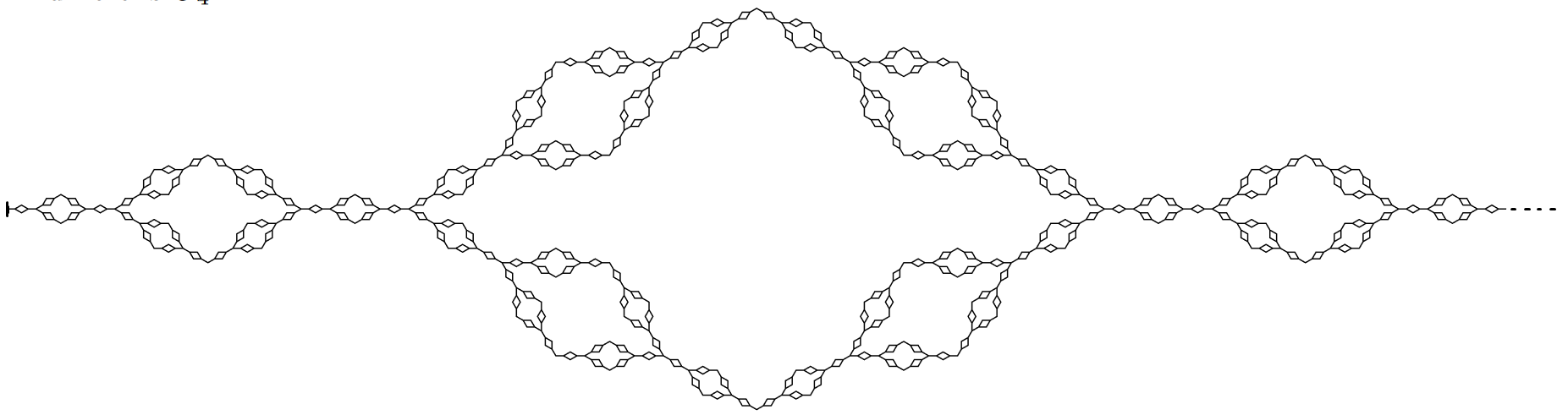}}; \node[circle, draw, fill=blue!50, inner sep=1.5pt, label=above:{$\origin$}] at (-.35\textwidth,0){};  
\end{tikzpicture}
\caption{Close up of a symmetrically self-similar graph with branching number 2 near its origin $\origin$. The generating cell is pictured in Figure~\ref{fig:cell}(a).}
\label{fig:diamond}
\end{figure}
In this work we focus on a particular family of fractal graphs,
introduced by Kr\"on called \emph{symmetrically self-similar
  graphs}~\cite{kron_green_2002}. These graphs are well studied,
\cite{malozemov_self_2003,bab_revisiting_2008,teufl_number_2011}, as
they are structured, and have many interesting properties.  We
give the precise definition in the next section, but simple examples
can provide some intuition. These graphs are generated by an iterative
process involving a finite \emph{cell graph} built from complete graphs
and have particular symmetry properties. The graph pictured in
Figure~\ref{fig:diamond} is generated from the cell in
Figure~\ref{fig:cell}(a), and the well known Sierpi\'nski graph is a
second example, generated from the cell in Figure~\ref{fig:cell}(b).
\begin{figure}
\centering 
\subfigure[$\theta=2$, $\mu=6$]%
{ 
\begin{tikzpicture}[scale=.75]\foreach \i/\x/\y in {
   1/0/0, 2/1/0, 3/2/0, 4/3/0, 
    5/1.5/.5, 6/1.5/-.5} 
{
    \node[circle, draw, inner sep=1.5pt] (\i) at (\x,\y) {};
}
    \node[circle, draw, fill=blue!50, inner sep=1.5pt, label=left:{$v_1$}] (1) at (0,0) {};
    \node[circle, draw, fill=blue!50, inner sep=1.5pt, label=right:{$v_2$}] (4) at (3,0) {};

\foreach \i/\j in {
1/2, 3/4, 
2/5, 5/3,  
2/6, 6/3}
{\draw[-] (\i) -- (\j);}

\end{tikzpicture}
}\quad
\subfigure[$\theta=3$, $\mu=3$]%
{
\begin{tikzpicture}[scale=.75]
    \node[circle, draw, fill=blue!50, inner sep=1.5pt, label=left:{$v_1$}] (v1) at (0,0) {};
    \node[circle, draw, fill=blue!50, inner sep=1.5pt, label=right:{$v_2$}] (v2) at (3,0) {};
    \node[circle, draw, fill=blue!50, inner sep=1.5pt, label=above:{$v_3$}] (v3) at (1.5,2.6) {};
    
    \node[circle, draw, inner sep=1.5pt] (v4) at (1.5, 0) {};
    \node[circle, draw, inner sep=1.5pt] (v5) at (0.75, 1.3) {};
    \node[circle, draw, inner sep=1.5pt] (v6) at (2.25, 1.3) {};
    
    \draw (v1) -- (v5) -- (v4) -- (v1);
    \draw (v3) -- (v5) -- (v6) -- (v3);
    \draw (v2) -- (v4) -- (v6) -- (v2);
\end{tikzpicture}
}\quad
\subfigure[$\theta=4$, $\mu=5$]%
{
\begin{tikzpicture}[scale=.75]
\foreach \i/\x/\y in {
    5/0/2, 6/1/3, 7/1/0, 8/0/1, 
    9/3/1, 10/2/0, 11/2/3, 12/3/2,
    13/1/2, 14/1/1, 15/2/1, 16/2/2} 
{
    \node[circle, draw, inner sep=1.5pt] (\i) at (\x,\y) {};
}
\foreach \i/\x/\y in {
    1/0/3, 2/0/0
   } 
{
    \node[circle, draw, fill=blue!50, inner sep=1.5pt, label=left:{$v_\i$}] (\i) at (\x,\y) {};
}
\foreach \i/\x/\y in {
   3/3/0, 4/3/3
   } 
{
    \node[circle, draw, fill=blue!50, inner sep=1.5pt, label=right:{$v_\i$}] (\i) at (\x,\y) {};
}
\foreach \i/\j in {
1/5, 1/6, 1/13, 
2/7, 2/8, 2/14, 
3/9, 3/10, 3/15, 
4/11, 4/12, 4/16, 
5/6, 5/13, 
6/13, 
7/8, 7/14, 
8/14, 
9/10, 9/15, 
10/15, 
11/12, 11/16, 
12/16, 
13/14, 13/15, 13/16, 
14/15, 14/16, 
15/16}
{\draw[-] (\i) -- (\j);}
\end{tikzpicture}
}

\caption{Three cell graphs, with their extremal vertices in blue.}
\label{fig:cell}
\end{figure}
Symmetrically self-similar graphs viewed as a discretization of
fractals appear in a variety of different domains including the
study of Brownian motion on fractals~\cite{lindstrom_brownian_1990}.

The Green's functions of these graphs are well
studied~\cite{grabner_functional_1997, kron_asymptotics_2004},
including the asymptotic behavior of the excursion
probabilities. Kr\"on and Teufl~\cite{kron_asymptotics_2004} show that
asymptotically, as $n$ tends towards infinity,
$p^{(n)}(\origin,\origin)$ tends to~$n^\beta w(n)$ where $\beta$ is
related to the geometry of the graph, and~$w(n)$ is some oscillating
function. In the case of the Sierpi\'nski graph, $\beta=\log{5}/\log{3}$. Because
$\beta$ is not rational, we can deduce, using the known structure of coefficients
asymptotics of solutions to differential
equations (see~\cite[Theorem 19.1 p. 111]{wasow_asymptotic_1987}), that the
Green's function is not D-finite. Indeed, in almost all examples that
have been studied, non-D-finiteness (and hence transcendence) of the
Green's function can be deduced immediately from the singular
expansion. The trivial fractal graphs consisting of a finite number of
lines radiating
from a single origin vertex (\emph{star graphs}; see
Figure~\ref{fig:star}) are the only known examples
of symmetrically self-similar graphs with algebraic Green's
functions. (Although, there are other families of fractal like graphs
with algebraic Green's functions.)

Recently, Di Vizio \emph{et al.}~\cite{di_vizio_inhomogeneous_2023}
showed that not only is the Green's function of the Sierpi\'nski graph not
D-finite, but it actually \emph{differentially  transcendental}!  In the
absence of analogous structure theorems for differentially algebraic functions, the proof uses the key result of Grabner and Woess~\cite{grabner_functional_1997} that the Green's function of the Sierpi\'nski graph satisfies an iterative equation of the form 
\begin{equation}
\label{eq:central}
G(z)= f(z)G(d(z)),
\end{equation}
where~$f$ and~$d$ are explicit rational functions and are themselves Green's functions of a finite graph related to the generating cell. By Di Vizio \emph{et al}., series solutions to such equations are either algebraic or are differentially transcendental. In the present work, we generalize this approach to the full class of symmetrically self-similar graphs.  This contributes to a growing effort to understand the combinatorial embodiment of differential algebraicity and differential transcendence. Walks in graphs have been a rich domain~\cite{dreyfus_nature_2018}, and Galoisian approaches on other examples have been successful~\cite{bostan_differential_2024}. Such questions have been around at least since Rubel's invitation~\cite{rubel_survey_1989} to the topic. 

%
 We conjecture the following, which directly connects differential transcendence to combinatorial structure.
\begin{conjecture}
\label{conj:main}
The Green's function of a symmetrically self-similar graph with bounded geometry is algebraic if, and only if, the graph is a star consisting of finitely many one-sided lines with exactly one origin vertex in common. Otherwise, it is differentially transcendental.
\end{conjecture}
We prove this conjecture for a subcase in the main theorem of this
paper, of which Figure~\ref{fig:diamond} is an example. The hypotheses
depend on the branching number $\theta$ of $X$, which is defined
below. Here we state our main theorem which appears in
Section~\ref{sec:main_theorem} as Theorem~\ref{thm:MAIN}. An example
of a star graph is in Figure~\ref{fig:star}.

\begin{nontheorem}[Main result] 
Let~$X$ be a symmetrically self-similar graph with bounded geometry, origin~$\origin$, and branching number~$\theta=2$. Either the graph is a star consisting of finitely many one-sided lines that coincide at $\origin$,
or the Green's function $G(\origin, \origin|z)$ of~$X$ is differentially transcendental over~$\mathbb{C}(z)$. 
\end{nontheorem}
\begin{figure}\center
  \begin{tikzpicture}[scale=.5]
    \draw[-] (5.1,0) -- (5.4,0){};
    \draw[-] (5.1,-1) -- (5.4,-1){};
    \draw[-] (5.1,-2) -- (5.4,-2){};
 \foreach \i/\x/\y in {
    0/0/0, 1/1/0, 2/2/0, 3/3/0, 4/4/0, 5/5/0,
    11/1/-1, 12/2/-1, 13/3/-1, 14/4/-1, 15/5/-1,
    21/1/-2, 22/2/-2, 23/3/-2, 24/4/-2, 25/5/-2} 
{ \node[circle, draw,fill=white, inner sep=1.5pt] (\i) at (\x,\y) {};
  }
\foreach \i/\j in {
  0/1, 1/2,2/3, 3/4, 4/5,
  0/11, 11/12, 12/13, 13/14, 14/15,
  0/21, 21/22, 22/23, 23/24, 24/25}
{\draw[-] (\i) -- (\j);}
   
    \node[circle, draw, fill=blue!50, inner sep=1.5pt,label=above:{$\origin$}] (1) at (0,0) {};
    \node[label=left:{$X=$}] (1) at (-1,0) {};
    \node[label=right:{$\cdots$}] (1) at (5.2,0) {};
    \node[label=right:{$\cdots$}] (1) at (5.2,-1) {};
    \node[label=right:{$\cdots$}] (1) at (5.2,-2) {};
\end{tikzpicture}
\caption{A close up of a star graph near the origin.}
\label{fig:star}
\end{figure}
The argument at the heart of the proof of our result has the potential
to extend to all branching numbers  with a little more information on
the singularities of $G$, as we shall discuss near the end. 

What is the potential impact of this result? Malozemov and
Teplyaev~\cite{malozemov_self-similarity_2003} showed that the
spectrum of the difference Laplacian on self-similar graphs consists
of the Julia set of a rational function and a (possibly empty) set of
isolated eigenvalues that accumulate to the Julia set.  Kr\"on and
Teufl conjectured that this Julia set is an interval if, and only if,
the graph is a star consisting of finitely many one-sided lines with
exactly one origin vertex in common, otherwise it is a Cantor
set. Our work connects the spectrum and the differential transcendency of the
Green's function. While we do not resolve the conjecture of Kr\"on and
Teufl, the generating function classification provides a new
perspective. 

Finally, as we alluded to earlier, some infinite Cayley graphs are also fractal, and their Green's functions satisfy functional equations of the type in Eq.~\eqref{eq:central}.     However, they are easily proved to be algebraic~\cite{bell_cogrowth_2023}. Can we characterize a wider class of fractal graphs with algebraic Green's functions?

\section{Symmetrically self-similar graphs}
\label{sec:self-sim-graphs}
There are two ways to view this family of graphs. Kr\"on describes a
blow-up construction applied to a cell in ~\cite{kron_green_2002}, and
this is further developed by and Kr\"on and
Teufl~\cite{kron_asymptotics_2004}. On the other hand, such a graph can
be defined as a fixed point, as was done by Malozemov and
Teplyaev~\cite{malozemov_pure_1995}.  We will consider both
definitions below, as they appeal to complementary intuitions.  

We start by recalling some basic graph concepts. Let  $X = (V(X), E(X))$ be an undirected, simple graph (no loops
nor multiple edges). Given a sub-graph $\cell$ of $X$, we define its
boundary in $X$, denoted $\theta\cell$, to be the set of all vertices
in $V(X)\setminus V(\cell)$ which are adjacent to some vertex in
$\cell$.
Define the graph $\hcell$ as the graph induced by the vertices of
$\cell$ and $\theta\cell$ inside $X$. 

\subsection{A fixed-point construction}
Let $B$ be a sub-graph of a simple graph~$X$.  We define the set
$\operatorname{Con}_X(B)$ to be the set of connected components in $X
\setminus B$. Associated to $B$, we define the graph $X_B$ with vertex
set $V(B)$ and an edge $\{x, y\}$ if there is a
component in $\operatorname{Con}_X(B)$ with both $x$ and $y$ in its
boundary in $X$. Roughly, they are both incident to a common component
created by deleting~$B$ from~$X$. We give a simple example below. 

\begin{definition}[self-similar with respect to a sub-graph, and a map]
A connected infinite graph~$X$ is said to be \emph{self-similar with respect to sub-graph}~$B$ and a map\[\psi: V(X) \to V(X_B)\] if the following hold:
\begin{enumerate}
	\item no two vertices in~$B$ are adjacent in~$X$;
	\item the intersection of the boundaries of two different components in $\operatorname{Con}_X(B)$ contains at most one vertex; and
	\item $\psi(X)=X_B$ is isomorphic to $X$.
        \end{enumerate}
\end{definition}
        
The components of $\operatorname{Con}_X(B)$ are called
\emph{cells}. Consider~ a map $\psi$ that replaces a cell~$\hcell$
in $X$ with a complete graph on the vertices of~$\theta \cell$. If
$\psi(X)$ is isomorphic to $X$, then we have self-similarity. 

\begin{example}
\label{ex:line}
Consider the infinite finite path graph~$X$ with every second vertex coloured
blue. 
\begin{center}
\begin{tikzpicture}[scale=1]\foreach \i/\x/\y in {
 2/1/0, 4/3/0, 6/5/0} 
{
    \node[circle, draw, fill=blue!50, inner sep=1.5pt] (\i) at (\x,\y) {};
  }
  \foreach \i/\x/\y in {
   3/2/0, 5/4/0}
{
    \node[circle, draw, inner sep=1.5pt] (\i) at (\x,\y) {};
  }
  
\node (1) at (.5,0) {};
\node (7) at (5.5, 0) {};
\foreach \i/\j in {
2/3, 3/4, 4/5, 5/6, 6/7}
{\draw[-] (\i) -- (\j);}
    \node[circle, draw, fill=blue!50, inner sep=1.5pt, label=above:{\hphantom{$\origin$}},label=below:{$\origin$}] (1) at (1,0) {};
    \node[label=left:{$X=$}] (1) at (.5,0) {};
    \node[label=right:{$\cdots$}] (1) at (5.5,0) {};
\end{tikzpicture}
\end{center}
Let $B\subset X$ be the set of blue vertices, then
$\operatorname{Con}_X(B)$ is the set of (isolated) white vertices. These are the
cells; remark the cells are all isomorphic. The closure of a cell includes the two adjacent blue
vertices. Each looks like: 
\begin{center}
\begin{tikzpicture}[scale=1]\foreach \i/\x/\y in {
   1/0/0, 2/1/0, 3/2/0} 
{
    \node[circle, draw, inner sep=1.5pt] (\i) at (\x,\y) {};
}
    \node[circle, draw, fill=blue!50, inner sep=1.5pt, label=left:{$\hcell=\quad v_1$}] (1) at (0,0) {};
    \node[circle, draw, fill=blue!50, inner sep=1.5pt, label=right:{$v_2$}] (3) at (2,0) {};

\foreach \i/\j in {
1/2, 2/3}
{\draw[-] (\i) -- (\j);}
\end{tikzpicture}
\end{center}
Let~$\psi$ be the map on~$X$ that replaces every instance of $\hcell$
in $X$ with a single edge between the blue vertices. The image of
$\psi$ is again an infinite finite path graph, and is hence in bijection with
$X$. Thus, we can say that $X$ is self-similar with respect to the set of blue vertices and $\psi$. 
\end{example}



\begin{definition}[Symmetrically self-similar graph]
A self-similar graph~$X$ is furthermore said to be \emph{symmetrically
self similar} if it is locally finite\footnote{all vertices have finite degree}, and it
additionally satisfies the following axioms:
\begin{enumerate}
	\item All cells are finite,  and for any pair of cells $C$ and $D$ in $C_X(F)$
		there exists a isomorphism on vertices from
		of $\hat C$ to $\hat D$ such that the image of $\theta C$ is $\theta D$;
	\item The automorphism group $\operatorname{Aut}(\hat C)$ of $\hat C$
		acts doubly transitively on $\theta C$,
		which means that it acts transitively on the set of ordered pairs
		 $\left\{ (x,y) \mid x,y \in \theta C, \; x \neq y \right\},$ 
		where $g((x,y))$ is defined as $(g(x),g(y))$ for any $g\in \operatorname{Aut}(\hcell)$.
\end{enumerate}
\end{definition}
We see that the finite path graph in the example is symmetrically
self-similar as the additional conditions are met. 

Under these conditions, an \emph{origin cell} is a cell $\cell$ such that $\psi(\theta C) \subset \cell$. A vertex fixed point of $\psi$ is called an \emph{origin vertex}.  If $X$ is self-similar with respect to $B$ and $\psi$ then by Kr\"on, (\cite[Theorem~1]{kron_green_2002})  $X$ has precisely one origin cell and no origin vertex or exactly one origin vertex, which we shall denote~$\origin$ when it is the case. 

As the cells are isomorphic, we can talk about the cell of a
symmetrically self-similar graph. Denote this cell by $\cell$ and its
boundary in $X$ by $\theta\cell$. The boundary vertices $\theta\cell$
are called the \emph{extremal vertices}. If there are $\theta$ of
them, then $\theta$ is the \emph{branching number} of $X$. If we write
$\hcell$ as \emph{the} cell associated to the graph, we are referring to the closure in $X$.

It can be shown that the edges of $\hcell$ can be partitioned
into~$\mu$ complete graphs on $\theta$ vertices, and that $\psi$ maps
instances of $\hcell$ to complete graphs on the $\theta$ extremal
vertices of $\hcell$ (i.e. $\psi(\hcell)$ is isomorphic to
$K_\theta$). By simple accounting we have
$\mu = \frac{2|E|}{\theta(\theta-1)}$ and so $\mu$ corresponds to the usual mass scaling factor of self-similar sets. 
Figure~\ref{fig:cell} has three examples of cell graphs  (indicating
the extremal vertices in blue), and showing the values of $\theta$ and
$\mu$.  If a self-similar graph is bipartite, its cell graph is
bipartite and $\theta=2$. 

\subsection{A blowup construction}
Krön and Teufl~\cite[Theorem 1]{kron_asymptotics_2004} describe the
reverse of this construction as an iterative process whose limit is a
symmetrically self-similar graph. The initial graph is $\hcell$, and
at each step the next graph is generated using an inverse action to
$\psi$: the $\mu$ copies of $K_\theta$ in $\hcell$ are ``blown up'',
and are replaced by a copy of the~$\hcell$.  The sequence of
substitutions converges to a unique graph $X$, which is fixed under $\psi$. 
The first few iteration that builds the Sierpi\'nski graph is pictured in Figure~\ref{fig:Serp}.

Remark that the graph with a finite number of copies of this graph
centered at $\origin$ is also a fixed point of $\psi$. That is,
given~$\hcell$, we can generate a countable family of symmetrically
self-similar graphs. However, remark that as copies are added, the
outward probability of the origin is reduced by the same factor, and
hence the Green's functions will all be the same. 
\begin{figure}
  \centering

$C_1=$ \begin{tikzpicture}[scale=.75]   
    \filldraw[fill=blue!20] (0,0) -- (0.75, 1.3) -- (1.5, 0) -- cycle;
    \filldraw[fill=blue!20] (1.5,2.6) -- (0.75, 1.3) -- (2.25, 1.3) -- cycle;
    \filldraw[fill=blue!20] (3,0) -- (1.5, 0) -- (2.25, 1.3) -- cycle;
    \node[circle, draw, fill=white, inner sep=1.5pt] (v1) at (0,0) {};
    \node[circle, draw, fill=white, inner sep=1.5pt] (v2) at (3,0) {};
    \node[circle, draw, fill=white, inner sep=1.5pt] (v3) at (1.5,2.6) {};
    \node[circle, draw, fill=white, inner sep=1.5pt] (v4) at (1.5, 0) {};
    \node[circle, draw, fill=white, inner sep=1.5pt] (v5) at (0.75, 1.3) {};
    \node[circle, draw, fill=white, inner sep=1.5pt] (v6) at (2.25, 1.3) {};
    

  \end{tikzpicture}%
  \hspace{2cm}%
$C_2=$ \begin{tikzpicture}[scale=.75]
    \node[circle, draw, inner sep=1.5pt] (v1) at (0,0) {};
    \node[circle, draw, inner sep=1.5pt] (v2) at (3,0) {};
    \node[circle, draw, inner sep=1.5pt] (v3) at (1.5,2.6) {};

    \node[circle, draw, inner sep=1.5pt] (v12) at (3,0) {};
    \node[circle, draw, inner sep=1.5pt] (v22) at (6,0) {};
    \node[circle, draw, inner sep=1.5pt] (v32) at (4.5,2.6) {};

   \node[circle, draw, inner sep=1.5pt] (v13) at (1.5,2.6) {};
    \node[circle, draw, inner sep=1.5pt] (v23) at (4.5,2.6) {};
    \node[circle, draw, inner sep=1.5pt] (v33) at (3,5.2) {};
    
    \node[circle, draw, inner sep=1.5pt] (v4) at (1.5, 0) {};
    \node[circle, draw, inner sep=1.5pt] (v5) at (0.75, 1.3) {};
    \node[circle, draw, inner sep=1.5pt] (v6) at (2.25, 1.3) {};

    \node[circle, draw, inner sep=1.5pt] (v42) at (4.5, 0) {};
    \node[circle, draw, inner sep=1.5pt] (v52) at (3.75, 1.3) {};
    \node[circle, draw, inner sep=1.5pt] (v62) at (5.25, 1.3) {};

    \node[circle, draw, inner sep=1.5pt] (v43) at (3, 2.6) {};
    \node[circle, draw, inner sep=1.5pt] (v53) at (2.25, 3.9) {};
    \node[circle, draw, inner sep=1.5pt] (v63) at (3.75, 3.9) {};
    
    \draw (v1) -- (v5) -- (v4) -- (v1);
    \draw (v3) -- (v5) -- (v6) -- (v3);
    \draw (v2) -- (v4) -- (v6) -- (v2);

    \draw (v12) -- (v52) -- (v42) -- (v12);
    \draw (v32) -- (v52) -- (v62) -- (v32);
    \draw (v22) -- (v42) -- (v62) -- (v22);

    \draw (v13) -- (v53) -- (v43) -- (v13);
    \draw (v33) -- (v53) -- (v63) -- (v33);
    \draw (v23) -- (v43) -- (v63) -- (v23);
\end{tikzpicture}

    \caption{ $C_1=\hcell$ and $C_2$. From $C_i$ to $C_{i+1}$ the
      shaded triangles are each replaced by a copy of $\hcell$. The limit graph in this case is the Sierpi\'nski graph.}
    \label{fig:Serp}
\end{figure}

%
%
\section{Green's functions}
\label{sec:green_functions}
Our work relies heavily on the existing literature on Green's functions. 
\subsection{Finite Graphs}
The study of Green's functions of finite graphs is a classic topic,
going back to at least the 1950s. We will recall some results as presented in~\cite{flajolet_analytic_2009}, as they are very helpful for our purpose.  Very roughly, in the case of a finite graph, we can express  Green's functions as quotient of two determinants of the probability matrix. Thus, in this case, the Green's function will be a rational function. 

\begin{prop}[Application of Proposition V.6 Flajolet and Sedgewick]
\label{thm:greenformula}  
Let~$W$ be a digraph on the vertex set~$\{v_1, \dots, v_k\}$ and let
$T$ be the $k\times k$ matrix of edge weights
$T[i,j]=\operatorname{weight}(v_i,v_j)$. The weight of a path is the
product of the weights of the edges in the path. 
The series $F_{i,j}(z)$ defined so the coefficient of $z^n$ is the sum of
the weights of all paths from $v_i$ to $v_j$ of length $n$ in $W$,  is the entry $[i,j]$ of the
matrix $(I-zT)^{-1}$. Otherwise stated,  
\begin{equation}
\label{eqn:greenformula}  
F_{i,j}(z)=\left( (1-zT)^{-1}\right)=(-1)^{i+j}\frac{\Delta^{j,i}(z)}{\Delta(z)}.
\end{equation}
Here $\Delta(z)=\det(I-zT)$ and $\Delta^{j,i}(z)$ is the determinant
of the minor of index $j,i$ of $I-zT$.
\end{prop}

Thus, in the finite case the Green's function is a rational function. We
can deduce further information about a singular expansion near its
dominant singularity.

\begin{lemma}[Lemma V.1  of \cite{flajolet_analytic_2009} (Iteration of irreducible
matrices).] Let the non-negative matrix~$T$ be irreducible and
aperiodic, with $\lambda_1$ its dominant eigenvalue. Then the residue
matrix such that
\[(I-zT )^{-1} = \frac{\Phi}{1-z\lambda_1} +O(1) \quad(z \rightarrow \lambda^{-1})\]
has entires given by $\phi_{i,j} = \frac{r_i\ell_j}{\langle r,
  \ell\rangle}$ where $r$ and $\ell$ are, respectively, right and left eigenvectors of $T$ corresponding to the
eigenvalue $\lambda_1$.
\end{lemma}

\begin{theorem}[Theorem V.7 of \cite{flajolet_analytic_2009} (Asymptotics of paths in finite graphs)] 
\label{thm:green_finite_graph}
Consider the matrix
\[ F(z)=\left(I-zT \right)^{-1}\]
where $T$ is a scalar non-negative matrix, in particular, the adjacency matrix of a
graph equipped with positive weights. Assume that $T$ is irreducible. Then all entries
$F_{i,j}(z)$ of $F(z)$ have the same radius of convergence~$\rho$, which can be defined in two
equivalent ways:
\begin{enumerate}
\item as $\rho=\lambda_1^{-1}$, where $\lambda_1$ is the largest positive eigenvalue of $T$;
\item as the smallest positive root of the determinantal equation: $\det(I-zT )=0$.
\end{enumerate}
Furthermore, the point $\rho=\lambda_1^{-1}$ is a simple pole of each $F_{i,j}(z)$.
If $T$ is irreducible and aperiodic, then $\rho=\lambda_1^{-1} $ is the unique dominant singularity
of each $F_{i,j}(z)$, and
\[[z^n]F_{i,j}(z)=\phi_{i,j} \lambda_1^n  +O(\Lambda^n), 0 \leq\Lambda<\lambda_1,\]
for computable constants $\phi_{i,j}>0$.
\end{theorem}

Given a finite graph with vertex set $\{v_1, \dots v_k\}$ if we set~$T$ to be~$P$, the matrix of transition probabilities as defined in Eq.~\eqref{eq:P},
then~$(T^n)_{ij}=p^{(n)}(v_i,v_j)$, and~$F_{ij}(z)=G(v_i, v_j|z)$. From this result, we can understand why
all Green's functions of finite graphs have the same radius of
convergence, specifically, the multiplicative inverse of the dominant eigenvalue of $T$.

Since our cells are connected, $T$ is irreducible, and it is only
periodic when the cell is bipartite. We have a good understanding of
this case, and principally, we can show (by considering $T^2$, which
will not be periodic) that in the value we are interested in:
\[[z^n]F_{1,1}(z)=\begin{cases}
                           \phi_{1,1} \lambda_1^n  +O(\Lambda^n), 0
                           \leq\Lambda<\lambda_1, & n\text{ even}\\
                           0, & \text{otherwise.}
                         \end{cases}\]
Recall that if $T$ is a stochastic matrix then 1 is a dominant eigenvalue of $T$.

\begin{example}[Diamond]
Consider the labelled graph and its the transition matrix is given in
  Figure~\ref{fig:graphandmatrix}.

\begin{figure}
\begin{minipage}{.45\textwidth}\center
  \begin{tikzpicture}[scale=1.2]
    \node[label={$\hcell=$}] at (-1.25, -0.25){};
    \node[circle, draw,  label=above left:{$w_{1}$},  inner sep=1.5pt] (w1)
    at (1,0) {};
    \node[circle, draw,  label=above right:{$w_{2}$},  inner sep=1.5pt] (w2)
    at (2,0) {};
    \node[circle, draw,  label=above:{$w_{3}$},  inner sep=1.5pt] (w3) at (1.5,.5) {};
    \node[circle, draw,  label=below:{$w_{4}$},  inner sep=1.5pt] (w4)
    at (1.5, -.5) {};
    \node[circle, draw, fill=blue!50, inner sep=1.5pt, label=left:{$v_1$}] (v1) at (0,0) {};
    \node[circle, draw, fill=blue!50, inner sep=1.5pt, label=right:{$v_2$}] (v2) at (3,0) {};

\foreach \i/\j in {
v1/w1, v2/w2, 
w1/w3, w1/w4,  
w2/w3, w2/w4}
{\draw[-] (\i) -- (\j);}
\end{tikzpicture}
\end{minipage}%
\begin{minipage}{.5\textwidth}\center
$P=\begin{array}{c|cccccc}
 &v_1 & v_2 & w_1& w_2& w_3& w_4\\  \hline    
v_1&0& 0& 1  & 0 & 0& 0\\
v_2&0 & 0 & 0&1 & 0& 0 \\
w_1&\frac13 & 0 & 0 &0&\frac13 &\frac13\\
w_2&0& \frac13&0&0& \frac13&\frac13\\
w_3&0&0&\frac12&\frac12&0&0\\
w_4&0&0&\frac12& \frac12&0&0
   \end{array}$
\end{minipage}   
\caption{A graph and its transition matrix}
\label{fig:graphandmatrix}
\end{figure}
This matrix has two  dominant eigenvalues 1 and -1, and we can see that 
\[G_{\hcell}(v_1, v_1\mid z)=F_{1,1}(z)=\frac{z^{4}-9 z^{2}+9}{3 \left(z^{2}-3\right)
    \left(z^{2}-1\right)}. \] The series starts $G(v_1, v_1\mid z)= 1+\frac{1}{3}
z^{2}+\frac{2}{9} z^{4}+\frac{5}{27} z^{6}+\frac{14}{81}
z^{8}+\mathrm{O}\! \left(z^{10}\right)$. For example, we see
that~$p^{(8)}(v_1, v_1)=\frac{14}{81}$.
Remark, the periodicity is evident as $G_{\hcell}(v_1, v_1\mid z)$ is
a function of $z^2$.
\end{example}

In the next section we apply these theorems to the cells of a
symmetrically self-similar graph, setting some edge weights to 0. 

\subsection{Symmetrically self-similar graphs}
Kr\"on and Teufl translate the iterative fractal nature of~$X$ into a
functional equation for the Green's function in one of the main results
of~\cite{kron_asymptotics_2004}, exploiting the iterative nature of
the blow up construction.  Define~$d$ to be the \emph{transition
  function}: the generating function of the probabilities that the
simple random walk on $\hcell$ starting at $\origin$ in $\theta\cell$
hits a vertex in $\theta\cell\setminus \{\origin\}$ for the first time
after exactly $n$ steps. It may return to the origin before hitting
this other cell.  Next, define $f$ to be the \emph{return function}:
the generating function of $f_n$, the probability that the random walk
on $\hcell$ starting at $\origin$ returns to $\origin$ after $n$
stages without hitting a vertex in $\theta\cell\setminus
\{\origin\}$. Since the start is considered the first visit, $f_0=1$,
and thus $f(0)=1$. We also use $r(z)$, the probability function of
first return which satisfies $f(z)=\frac{1}{1-r(z)}$.

Because~$\hcell$ is finite, $f$ and $d$ are expressed in terms of Green's
functions of a finite graph. They are well understood: they are
rational functions that can be computed as the determinant of a
specific matrix. We first give the functional equation, and then
consider singular expansions of~$G$,~$d$ and~$f$.

Given the closed cell~$\hcell$, and let $P_{\hcell}$ be the
probability transition matrix associated to it.
Set $P_f$ to be the matrix $P_{\hcell}$ with the modification
$P_{f}[i, j]:=0$ for $1<j\leq \theta$, $i>\theta$. Then 
$f(z)=\left[(I-zP_f)\right]^{-1}[1,1]$.  Set $P_d$ to be the matrix $P_{\hcell}$ with the modification
$P_{d}[i,j]:=0$ for $1<i\leq \theta$, $j>\theta$. Then $d(z)=\sum_{j=2}^\theta\left[(1-zP_d)\right]^{-1}[1,j]$.

\begin{lemma}[Kr\"on and Teufl~\cite{kron_asymptotics_2004}, Lemma 3]
\label{lem:GSS}
Let~$G(z)$ be the Green's function for a symmetrically self-similar graph with origin $\origin$. Then, the following equation holds for all $z \in U(0,1)$:
\begin{equation}
    G(z) = f(z)G(d(z)).
\end{equation}
\end{lemma}
This functional equation was first established and exploited in the
case of the Sierpi\'nski graph by Grabner and Woess
\cite{grabner_functional_1997} before being generalized for all
self-similar graphs by Krön \cite{kron_green_2002}.


\begin{example}[Star graphs]
Consider at cell that is a path on three vertices, as in Example 1. It is
straightforward to show that $f(z)=\frac{2}{2-z^2}$ and
$d(z)=\frac{z^2}{2-z^2}$, and that any star graph thus has $G(z)= \frac{1}{\sqrt{1-z^2}}$, as
it satisfies $G(z)=f(z)G(d(z))$. Remark this is an algebraic function,
a solution $y$ to $y^2(1-z^2)-1=0$. We further deduce $p^{(n)}(\origin, \origin) =\binom{2n}{n}2^{-n}$,
for even $n$, 0 otherwise.
\end{example}

\begin{example}[Diamond example continued]
  We can define $f$ and $d$ for the diamond example using the
  following matrices of graph weights:
  \begin{equation}
    P_f=\left[\begin{array}{cccccc}
0& 0& 1  & 0 & 0& 0\\
0 & 0 & 0&1 & 0& 0 \\
\frac13 & 0 & 0 &0&\frac13 &\frac13\\
0& 0&0&0& \frac13&\frac13\\
0&0&\frac12&\frac12&0&0\\
0&0&\frac12& \frac12&0&0
\end{array}\right]
        \quad
 P_d=\left[\begin{array}{cccccc}
0& 0& 1  & 0 & 0& 0\\
0 & 0 & 0&0 & 0& 0 \\
\frac13 & 0 & 0 &0&\frac13 &\frac13\\
0& \frac13&0&0& \frac13&\frac13\\
0&0&\frac12&\frac12&0&0\\
0&0&\frac12& \frac12&0&0
\end{array}  \right].
\end{equation}
It follows that
\begin{eqnarray*}
 f(z) &= (1-zP_f)^{-1}[1,1]
=-\frac{3 \left(2 z^{2}-3\right)}{z^{4}-9 z^{2}+9}=1+\frac{1}{3}
z^{2}+\frac{2}{9} z^{4}+\mathrm{O}\! \left(z^{10}\right)\\
\text{and}\quad d(z)&= (1-zP_d)^{-1}[1,2]=\frac{z^{4}}{z^{4}-9 z^{2}+9}
                             =\frac{1}{9} z^{4}+\frac{1}{9}
                             z^{6}+\frac{8}{81} z^{8}+\mathrm{O}\!
                             \left(z^{10}\right).
\end{eqnarray*}
The dominant eigenvalue of both $P_f$ and $P_d$ are
$\frac{1+\sqrt{5}}{2\sqrt{3}}$ and hence $\rho_f=\rho_d =
\frac{2\sqrt{3}}{1 + \sqrt{5}} \approx 1.070$. Using these, we can
determine the initial series expansion of the Green's function,
\begin{equation}
                             G(\origin, \origin|z) \equiv G(z)=
                             f(z)\,f(d(z))\,f(d(d(z)))\dots =1+\frac{1}{3}
                             z^{2}+\frac{2}{9} z^{4}+\frac{5}{27}
                             z^{6}+O(z^7).
\end{equation}
\end{example}  

We inspect the components of this equation more closely.
The Green's function $G(u,v|z)$ of a finite graph with transition
probability matrix $P$ is the corresponding $[u,v]$ entry of
$(I-zP)^{-1}$. Analytically, we know a fair amount about $G$ in this
case.  Remark the result does not apply directly to bipartite graphs,
as $P$ fails the condition of aperiodicity, but there are workarounds
using $P^2$, for example. A key result is the singular expansion of $f$ (and $d$) is of the form $\kappa(1-z/\rho)^{-1} + O(1)$,  as $z\rightarrow \rho$ for a real, positive $\kappa$. Furthermore, as $\hcell$ is finite, $d$ and $f$ are rational. 

A graph with bounded geometry here refers to a graph where the number
of neighbours each vertex has is finite and has a uniform upper bound.
In the case of symmetrically self-similar graphs with bounded geometries Kr\"on and Teufl exploit their functional equation to determine the singular expansion of~$G(\origin,\origin|z)$ near~$z=1$. Here $\tau=d'(1)$ and $\alpha=f(1)$.
\begin{theorem}[Kr\"on and Teufl~\cite{kron_asymptotics_2004}, Theorem 5]
\label{thm:KT}
Let~$X$ be a symmetrically self-similar graph with bounded geometry and origin vertex~$\origin$.  Then there exists a 1-periodic, holomorphic function $\omega$ on some horizontal  strip around the real axis such that the Green's
function $G$ has the local singular expansion
\begin{equation}
\label{eq:KT}
\forall |z|<1, \quad G(z) = (z-1)^{\eta} \left( \omega \left( \frac{1-z}{\tau} \right) + o(z-1) \right) 
\end{equation}
 with  $\eta = \frac{\log \mu}{\log \tau}-1$.
\end{theorem}
Recall that~$\mu$ is the number of cliques that form~$\hcell$. When $\theta=2$ this is the number of edges in $\hcell$. Moreover, even though it is not directly stated, it follows from a quick analysis of the asymptotics of~$G$ at~$1$~that $\eta = -\frac{\log \alpha}{\log \tau}$. 

\begin{example}[Diamond example continued]
  Since we have exact expressions for $f$ and $d$ we compute $\eta$:
  \[\eta = \frac{\log\mu}{\log\tau}-1=\frac{\log 6}{\log 18}-1=-\frac{\log
 3}{\log18}=\frac{f(1)}{d'(1)}\approx -0.3800.    \]
\end{example}

The domain of analyticity of $G$ is governed by the dynamic behaviour
of $f$ and $d$, according to Kr\"on~\cite{kron_green_2002} and Kr\"on
and Teufl~\cite{kron_asymptotics_2004}. If~$\omega$ is not constant, a
singularity analysis shows that $G$ is not algebraic. Kr\"on and Teufl conjecture that: if $\omega$ is constant, then
$\hcell$ is a path graph.

More is true, thanks to a recent result of Di Vizio, Fernandes and Mishna~\cite{di_vizio_inhomogeneous_2023}. 
\begin{theorem}[Di Vizio, Fernandes, and Mishna~\cite{di_vizio_inhomogeneous_2023}, Special case of Theorem C]\label{thm:DFM}
Suppose $y(z)\in\mathbb{C}[[z]]$ is a series that satisfies the
following iterative equation with some $a, b\in\mathbb{C}(z)$:
\[y(z)=a(z)y(b(z)).\]
    If $b(z)\in \mathbb{C}(z)$ satisfies the following conditions : 
$b(0) = 0$;  $b'(0) \in \{ 0,1, \emph{roots of unity}\}$; and
 no iteration of $b(z)$ is equal to the identity, then either there is some $ N \in \mathbb{N}^*$, such that $y(z)^N$ is rational,  or $y$ is differentially transcendental over $\mathbb{C}(z)$. 
\end{theorem}
Putting these results together we have our workhorse result. 
\begin{corollary}
\label{cor:DFM}
Let $G(z)$ be the Green's function of a symmetrically self-similar
graph with origin $\origin$. Then, $G(z)$ is either algebraic or
differentially transcendental over $\mathbb{C}(z)$.  If $G$ is
algebraic, then there exists a minimal $N$ such that $G^N=P/Q$
with~$P, Q\in\mathbb{C}[z]$, $P$ and $Q$ co-prime and $Q$ monic. In
this case, the equality
\begin{equation}
  \label{eqn:func-eqn}
  f(z)^N\frac{P(d(z))}{Q(d(z))}=\frac{P(z)}{Q(z)}
\end{equation}
is defined and true for all complex~$z$ except possibly a finite set
of points. 
\end{corollary}

\begin{proof}
First we recall that there is a unique analytic continuation of any Green's function $G(u,v|z)$ to some domain beyond the radius of convergence. 
If $G$ is algebraic, using the polynomial equation it satisfies we can verify that Eq.~\eqref{eq:central} holds in that domain. 

Next, under these circumstances, we verify that the hypotheses of
Theorem~\ref{thm:DFM} hold. We set $b$ to $d$, $f$ to $a$. These are
both rational since the cells are finite. We then set $y$ to $G$. By
definition, two vertices on the boundary of the cell have distance at
least two away, so $0=d_0=d(0)$ and $0=d_1=d'(0)$. No iteration of
$d(z)$ is the identity, since each iteration increases the number of
initial terms equal to 0. The rest follows from Theorem~\ref{thm:DFM}
since the hypotheses hold.

As to the validity of Eq.~\eqref{eqn:func-eqn}, this follows from the
fact that it can be reduced to a polynomial equation that is valid in
the entire open unit disc. The exceptions are contained in the set of
poles of $f$ and $d$, the zeroes of $Q$ and any $z$ such that $d(z)$
is a zero of $Q$. 
\end{proof}

\subsection{A few technical results on $f$ and $d$}
The derivative $d'(1)$ is called $\tau$ and is the expected number of
steps for a walk on $\hcell$ that starts at
$\origin$ to hit another external vertex. By definition, external
vertices are not adjacent, so $d'(1)>2$. 
Since $d_n$, ($f_n$ and $r_n$) are all probabilities, and hence
non-negative numbers, by Pringsheim's theorem the radius of
convergence $\rho_d$ (respectively $ \rho_f$ and $\rho_r$) of $d$,
($f$ and $r$) are singularities. Now, as $d$ is a probability
generating function it follows that $d(1)=1$ and thus $\rho_d>1$ since $d$ converges at~1. 
As these are Green's functions, we can deduce key facts about these elements that we use in subsequent proofs. The proofs are a mix of basic facts on Green's functions, and series analysis.
\begin{lemma} 
\label{lem:technical}
For $\rho_d, \rho_f$ and $\rho_r$ as defined above the following are true: 
\begin{enumerate}
\item $\rho_d=\rho_f$ \label{lem:rhodisrhof}
\item $\rho_d$ and $\rho_f$ are each poles of order 1; 
\item $\rho_f < \rho_r$; \label{lem:rhoflessthanrhor}
\item $\forall z\in (1, \rho_d), \, d'(z) > 1$ and therefore $d(z)>z$;
\item for $z \in [0, \rho_r]$ if $f(z) = \infty$, then $z=\rho_f$. \label{lem:rhoflessthanrhor}
\end{enumerate}
\end{lemma}

\begin{proof}[Proof of Lemma~\ref{lem:technical}]\mbox{}
\begin{enumerate}
\item 
We apply Eq.~\eqref{eqn:greenformula} to $P_f$ and $P_d$ to show
that the radii of convergence of $f$ and $d$ are the same. These two
matrices differ in the following ways -- The columns of $P_f$
associated to $v_i$, $i>1$ are zero since there are no edges to these
vertices. The rows of $P_d$ associated these $v_i$ are zero. The
determinants of $I-zP_f$ and $I-zP_d$ can be seen to be the same by
iteratively taking cofactor expansions recursively along the rows
associated to $v_i$  and the columns associated to $v_i$ respectively. In both cases,
the only non zero element is the 1 on the diagonal, and up to a sign
the value of the determinant is equal to the determinant of the matrix with that row and column
removed. Iterating this process we can see that as $P_d$ and $P_f$
with all columns and rows associated to the $v_i$ are removed give the
same matrix, both determinant computations will eventually reduce to the same value. 
\item This is a direct application of~\cite{flajolet_analytic_2009}, Theorem V.7.\emph{iii} and V.7.\emph{i}.
\item 
Since they are sums of non-negative terms, $f$ and $r$ are continuous increasing functions along the interval from 0 to their radius of  convergence. The inequality is a consequence of the relation   $f(z)=\frac{1}{1-r(z)}$: $f(z)$ is singular at the minimum of
$\rho_r$ and a value where $r(z)=1$. Since~$r(0)=0$, and $r(z)$ is increasing, the value of $z$ where it is 1 is strictly less than $\rho_r$.
\item We note that for $z$ in the radius of convergence, 
$d''(z) = \displaystyle \sum_{n \leq 2} n(n-1)d_n z^{n-2} > 0$, therefore, $d'$ is increasing. And $d'(1) = \tau \geq 2$.   
\item 
This also follows directly from  $f(z)=\frac{1}{1-r(z)}$, given that $r(\rho_f)=1$, and $r(z)$ is increasing up to its radius of convergence. 
\end{enumerate}
\end{proof}

\section{Proof the Main Theorem}
\label{sec:main_theorem}
Throughout this section, we assume~$X$ is a symmetrically self-similar
graph with bounded geometry, origin $\origin$ and branching number
$\theta=2$. If ever we additionally assume that
$G(z)=G(\origin, \origin|z)$ is algebraic, then we write $G^N=P/Q$ as
in Corollary~\ref{cor:DFM}. Recall the statement of the Main Theorem:
\begin{theorem}
\label{thm:MAIN}
Let $X$ be a symmetrically self-similar graph with bounded geometry,
origin $\origin$ and branching number $\theta=2$. Either the graph is
a star consisting of finitely many one-sided lines coinciding at
$\origin$, or the Green's function $G(\origin, \origin|z)$ is differentially transcendental over $\mathbb{C}(z)$. 
\end{theorem}
\begin{proof}[Outline of Proof]
First suppose $X$ is bipartite. We show in Theorem~\ref{lem:bipartite}
below that if $G$ is algebraic, then $\hcell$ is a finite path
graph. By~\cite[Lemma 2]{kron_green_2002}, $X$ is a star centered at
an origin vertex precisely when~$\hcell$ is a path. Otherwise, $G$ is
differentially transcendental.

Otherwise, if $X$ is not bipartite, we show in Theorem~\ref{lem:non-bipartite}, that $G$ is differentially
transcendental. Putting all of these components together, the result is established. 
\end{proof}

We can reframe this result as a case of the conjecture of Kr\"on and Teufl.
\begin{corollary}\label{cor:KT}
  Let~$X$ be a symmetrically self-similar graph with bounded geometry,
  origin~$\origin$, and branching number~$\theta=2$. Suppose that
  $\omega$ in Theorem~\ref{thm:KT} is constant, and $G$ is
  algebraic. Then, the graph is a star consisting of finitely many
  one-sided lines coinciding at $\origin$.
\end{corollary} 
%
\subsection{The structure of an algebraic $G(z)$}
\label{sec:algebraic}
We start by developing the properties of~$G$ in the case that it is
algebraic. 
\begin{lemma}
\label{lem:no_positive_poles}
If $G(z)$ is algebraic, then the rational function $G^N$ has only real
zeroes and poles, and no real poles nor zeroes greater than 1. 
\end{lemma}
\begin{proof}
If~$G$ is algebraic, then $\omega$ in Theorem~\ref{thm:KT} is constant
and by Kr\"on~\cite[Theorem 7]{kron_green_2002} the singularities of
$G(z)$ are real, and contained in the set $(-\infty, -1]\cup[1,
\infty)$.  Writing $G^N=P/Q$ we have that the singularities of $G$ are
contained in the the roots of $P$ and $Q$, as these are the poles and zeroes of $G^N$.

Towards a contradiction, assume there is some real $z_0$, with $z_0>1$
that is either a pole or a zero of $G^N$. Since $G^N$ has only a
finite number of zeroes and poles, we can assume that $z_0$ is the
smallest one strictly greater than $1$. Since $\rho_d$ is a pole of
$d(z)$,  $d((1,\rho_d)) = (1,\infty)$, and so there exists a $z' \in
(1,\rho_d)$ such that $d(z') = z_0$. Similarly, $f$ analytic and
defined on the interval $[1, \rho_d)$. By Lemma~\ref{lem:technical}, $\rho_d=\rho_f$, so $z'$
is not a pole of $f$. Since $f(0)=f_0=1$ (by convention, the length of
a walk whose first contact with the origin for a walk that starts at
the origin is 0) and $f$ is increasing on this interval so $z'$ is not a zero of $f$ either.

From Eq.~\eqref{eqn:func-eqn}, $\frac{G(z)^N}{f(z)^N} = G(d(z))^N$, and so
\begin{equation}
\lim_{z\rightarrow z_0}G(z)^N
=  \lim_{z\rightarrow  z'}G(d(z))^N
=\lim_{z\rightarrow z'} \frac{G(z)^N}{f(z)^N}.
\end{equation}
Hence, as $\lim_{z\rightarrow z'} f(z)=f(z')$ is finite and nonzero, if $z_0$ is a zero or pole of $G^N$ greater than~1, then~$z'$ is also a zero or pole, but is closer to 1, a contradiction. This contradiction establishes the result.   
\end{proof}

\begin{lemma}
\label{lem:deg}
In this context, if $G$ is algebraic, with $G^N=P/Q$, then~$\deg Q\geq N$. 
\end{lemma}
\begin{proof}
We start with Eq.~\eqref{eqn:func-eqn} and consider expansions of $G$, $f$, and $d$ around $\rho_f=\rho_d$. Now, $G(z)$ is not singular at $\rho_f$, since $\rho_f>1$, nor is $G$ zero at $\rho_f$, by Lemma~\ref{lem:no_positive_poles}.
Here $\kappa_f$ and $\kappa_d$ are positive real constants by~\cite[Theorem 7.]{flajolet_analytic_2009}. We compute the following limits in~$\mathbb{R}$:
\begin{eqnarray*}
G(\rho_f)^N &=&\lim_{z \rightarrow \rho_f} G(z)^N \\
&=&\lim_{z \rightarrow \rho_f} f(z)^N \frac{P(d(z))}{Q(d(z))}.\\
&=&\lim_{z \rightarrow
    \rho_f}\left(\frac{\kappa_f^N}{(1-z/\rho_f)^{N}} +
    O(z-\rho_f)^{-N+1}\right) \frac{P(\kappa_d(1-z/\rho_f)^{-1} +
    O(1))}{Q(\kappa_d(1-z/\rho_f)^{-1} + O(1))}\quad\text{by Lemma~\ref{lem:technical},}\\
&=&\lim_{u\rightarrow \infty} (\kappa_f^N u^{N} + O(u)^{N-1})
    \frac{P(\kappa_d u + O(1))}{Q(\kappa_d u + O(1))}\quad\text{with }u=(1-z/\rho_f)^{-1}\\
&=&\lim_{u\rightarrow \infty} \kappa_f^N u^{N}  \frac{P(\kappa_d u)}{Q(\kappa_d u)}.
\end{eqnarray*}
Remark as $z$ approaches $\rho_f$, $u$~is arbitrarily large. The degree of $Q$ as a polynomial in~$z$ must be larger than~$N$ since the value of the last limit is non-zero.
\end{proof}

\subsection{Results when $\theta=2$}
\label{sec:thetaistwo}
\begin{theorem}
\label{thm:alpha_mu}
If $\theta=2$ then  $\alpha \leq \mu$ with equality if and only if $\hcell$ is a finite path graph. 
\end{theorem}
The proof is a little technical, and is delayed to the next
section. Recall that $\mu$ is the number of edges in $\hcell$, and
$\alpha$ is the expected number of visits to $v_1$ before a visit to
$v_2$. The proof but relies on a straightforward intuition that
among all finite cell graphs $\hcell$ where $v_1$ and $v_2$ are a
fixed distance apart, the average number of times a random walk
starting at $v_1$ hits $v_1$ before $v_2$ (i.e. $f(1)$) is minimized
when the graph is a line.

\begin{corollary}\label{lem:eta}
    If $\theta=2$, then $\eta$ in the development in Eq.~\eqref{eq:KT} satisfies \[\eta \geq -\dfrac{1}{2},\] with equality if, and only if, $\hcell$ is a finite path graph. 
\end{corollary}
\begin{proof}
Starting from the definition of $\eta$ we compute the following bound:
\begin{eqnarray*}
\eta =    -\frac{\log \alpha}{\log \tau}
  &= &     \frac{\log \mu}{\log \tau} -1 \\
  &\geq& \frac{\log \alpha}{\log \tau} - 1\quad\text{since }\quad\mu\geq \alpha\\
  &=&     -\eta -1.
\end{eqnarray*}  
Solving for $\eta$ gives $\eta \geq -\frac{1}{2}.$
\end{proof}

\subsection{Case 1: $X$ is bipartite}
\label{sec:bipartite}
First we consider the case of bipartite graphs. In this case it must be that~$\theta=2$, since any complete graph on three or more vertices is not bipartite.  The main result of this subsection is the following: 
\begin{theorem}
\label{lem:bipartite}
Let $X$ be a symmetrically self-similar graph with bounded geometry,
origin $\origin$ and branching number $\theta=2$. If, additionally $X$ is bipartite, and $G$ is algebraic, then $\hcell$ is a path.  
\end{theorem} 
\begin{proof}
By Corollary~\ref{cor:DFM}, either $G$ is differentially
transcendental, or it is algebraic. Assume that $G$ is algebraic. By
Lemma~\ref{lem:no_positive_poles} the only positive singularity of $G$ is at~1, and by
Theorem~\ref{thm:KT} the exponent at the singular expansion is~$\eta$. The
singularities of $G$ are contained in the roots of $P$ and $Q$. Since
they are co-prime $P$ and $Q$ have no common roots, we deduce that
both $1$ and $-1$ are roots of order $-nN$ of $Q$. Up to a constant,
$Q(z)$ is $(1-z)^{-\eta N}(1+z)^{-\eta N}$.
Now, $\deg Q = -2\eta N$ and $\deg Q \geq N$ from
Lemma~\ref{lem:deg}. Putting the two relations together we have $\eta
\leq -1/2.$ This, coupled with Corollary~\ref{lem:eta} suggests that
$\eta=-1/2$ precisely.  By Corollary~\ref{lem:eta}~this equality holds
if and only if $\hcell$ is a path. 
\end{proof}  

\subsection{Case 2: $X$ is not bipartite}
\label{sec:nonbipartite}
Next we consider the case that $X$ is not bipartite.
\begin{theorem}\label{lem:non-bipartite}
Let $X$ be a non-bipartite symmetrically self-similar graph with origin $\origin$, bounded geometry and branching number $\theta=2$. Then the Green's function $G(z)$ is differentially transcendental over $\mathbb{C}(z)$. 
\end{theorem}
\begin{proof}
We assume $G$ is algebraic and derive a contradiction, and thus $G$
must be differentially transcendental.
If $X$ is not bipartite, then $z=1$ is the only singularity on the
boundary of the unit disc by~\cite[page
13]{kron_asymptotics_2004}. Since $G$ is algebraic, then $\omega$ in
Theorem~\ref{thm:KT} is constant, and thus the singularities of $G$
are contained in the set $(-\infty, -1)\cup[1, \infty)$ but we have
already shown that are no singularities in $(1, \infty)$ by Lemma~\ref{lem:no_positive_poles}. 
We next show that there are also no singularities in $(-\infty, -1)$.
Since $\rho_d$ is a pole of order $1$ of $d$, and $d$ is increasing,
we can deduce that $\lim_{z \rightarrow \rho_d^+} d(z) =
-\infty$. Using a limit argument applied to Eq.~\eqref{eqn:func-eqn}, We can
show $f$ has neither zero nor pole in the interval~$(\rho_f, \rho_r)$,
and $d$ also has no pole in the interval $(\rho_f, \rho_r)$.  
Therefore, by intermediate value theorem, for any $z_0 \in (-\infty,
-1)$, there is some $z' \in (\rho_f, \rho_r)$ such that $d(z') =
z_0$. Since $z' >\rho_f$, $z'$ is a positive real number greater than
1,  therefore, $z'$ is neither a pole nor a zero of $G$.

We substitute this into Eq.~\eqref{eqn:func-eqn}
\[ G^N(z') = \frac{P(z')}{Q(z')}=
  f(z')^N\frac{P(d(z'))}{Q(d(z'))}=f(z')^N\frac{P(z_0)}{Q(z_0)} =
  f(z')G(z_0).\]
We can conclude that $G(z_0)$ is finite, and nonzero.

Thus, given the singular expansion of $G$ at 1, we have that $Q$ has a
single zero and it is at 1, and is of order $-N\nu$. Combining prior
results:
\[
  N\leq \deg Q = -N\eta \implies \eta \leq -1.
\]
This directly contradicts that $\eta\geq -\frac{1}{2}$. Thus, $G$
cannot be algebraic.
\end{proof}

Remark that for $\theta>2$, it suffices to find a similarly
incompatible lower bound on $\eta$. This same proof should apply. 
\section{Proof of Theorem~\ref{thm:alpha_mu}}
To finish,  we now prove Theorem~\ref{thm:alpha_mu}  which states that
if $\theta = 2$, then $\mu \geq \alpha$ with equality if and only if
$\cell$ is a finite path graph. In fact, we conjecture that the
(strict) inequality holds for all cell graphs, even when $\theta > 2$. 
 
\subsection{Harmonic functions on graphs}
Let $\hcell$ be a cell of a symmetrically self-similar graph with $\theta = 2$.  For $v \in V(\hcell)$,  let $H(v)$ be the probability that a random walk starting at $v$ will reach $v_1$ before $v_2$. 
Hence $H(v_1) = 1$ and $H(v_2) = 0$. 
By decomposing by the first step of the random walk, we have 
\[H(v) = \dfrac{1}{\deg u} \displaystyle \sum_{u \in N(v)} H(u).\]
Then $F$ is a harmonic function of the graph $C$ with boundary conditions given by the values on $v_1$ and $v_2$. It is a classic result  that there is a unique harmonic function on a graph with fixed boundary conditions. 
Since $\theta = 2$, $v_1$ has only one neighbour,  call it $w$. Then $\alpha-1$ is the average time that a random walk starting at $v_1$  makes the step $w \rightarrow v_1$. After that step, the only possible step is the step $v_1 \rightarrow w$. 
A random walk on $\hcell$ starting at $v_1$ can has a unique decomposition\footnote{See Flajolet and Sedgewick~\cite{flajolet_analytic_2009} for a detailed treatise on combinatorial grammars that formalize such decompositions. } on the steps from $v_1$ to $w$ (which we denote $v_1 \rightarrow w$):
\[\mathcal{W}_1 = (v_1 \rightarrow w) \times \mathbf{SEQ}(\mathcal{W}_2 \times (v_1 \rightarrow w)) \times \mathcal{W}_3.\]
Here 
\begin{itemize}
    \item $\mathcal{W}_1$ is the set of walks starting at $v_1$ that stop after reaching $v_2$, but can return to~$v_1$ multiple times. 
    \item $\mathcal{W}_2$ is the set of walks starting at $w$ that never reach $v_2$ and stop after reaching $v_1$. 
    \item $\mathcal{W}_3$ is the set of walks starting at $w$ that never reach $v_1$ and stop after reaching $v_2$. 
\end{itemize}
We denote by~$\alpha$ the expected number of visits to $v_1$ in a walk in $\mathcal{W}_1$ (including the start of the walk). From the decomposition we see that this is one more than the expected number of $\mathcal{W}_2$ sub-walks, as each one ends in a unique return to~$v_1$.

Let $\mathcal{W}_4$ be the set of walks starting at $w$ that stop the first time they reach either $v_1$ or $v_2$. We can decompose this set by endpoint ($v_1$ and $v_2$ respectively) into disjoint sets:
\[\mathcal{W}_4 = \mathcal{W}_2 + \mathcal{W}_3.\]

The probability that a walk in $\mathcal{W}_4$ is also in $\mathcal{W}_2$ is precisely $H(w)$.  Therefore, the number of $\mathcal{W}_2$ in the decomposition of walk in $\mathcal{W}_1$ follows a geometric law of parameter $H(w)$.  Hence, $\alpha-1 = \dfrac{1}{1-H(w)}$. Finding an upper bound of $\alpha$ is therefore equivalent to finding an upper bound of $H(w)$. 

\subsection{The finite path graph}
We can determine $\alpha$ explicitly for a path $L$ on $\mu$ edges, as the harmonic function $H_L(v)$ is straightforward to compute. Remark $d(v_2,v_1)=\mu$ in this case.
Classically, 
\begin{equation}\label{eq:F_line}
H_L(v) = \frac{d(v_2,v)}{d(v_2,v_1)}.
\end{equation}
Remark, $H_L(v_1) = 1$, $H_L(v_2) = 0$ and for the vertex $v$ at distance $d$ from $v_2$, $d\in \{1, \cdots, d(v_1,v_2)-1 \}$, 
\[H_L(v) = \frac12\left(\frac{d-1}{d(v_1, v_2)}+\frac{d+1}{d(v_1, v_2)}\right)=\frac{d}{d(v_1, v_2)}.\] 
Therefore, in the case of a path, \[\alpha = 1+\frac{1}{1-H_L(w)}=1+ \frac{1}{1 - \frac{d(v_1,v_2) -1}{d(v_1,v_2)}} = d(v_1,v_2) = \mu.\]

\subsection{General case}
Consider a cell graph $\hcell$. We compare values of $H$, the harmonic function defined on $\hcell$ and boundary values $H(v_1)=1, H(v_2)=0$ to~$H_L$ as defined in Eq.~\eqref{eq:F_line}. 

First consider the following lemma:
\begin{lemma}
    Let $\hcell$ be a symmetric graph with $\theta = 2$ and $F$ be the harmonic function defined on $\hcell$ the value on the boundary: $H(v_1) = 1$ and $H(v_2)=0$. Let $u \in V(\cell)$ be a vertex such that there is a connected component $\tilde{C}$ of, $C \setminus \{u\}$ which contains both $v_1$ and $v_2$. Let $\tilde{H}$ be the harmonic function defined on $\tilde{C}$ by the same boundaries as for $F$. Then, 
    \[\forall x \in V(\tilde{C}), \quad H(u) \leq H(x) \Longrightarrow H(x) \leq \tilde{F}(x).\]
\end{lemma}
That is, if the probability of a reaching $v_1$ before $v_2$ is greater starting at $x$ than at $u$, this probability is bounded above the analogous probability in the graph where $u$ is removed, provided that $v_1, v_2$ and $x$ are all in the same component after the removal.
\begin{proof}
    The walks in $\hcell$ starting at $x$ absorbed by a boundary can be decomposed into three disjoint sets, according to which point $v_1, v_2, u$ it will reach first. We assign the probabilities,  $a,b,1-a-b$, to be the respective probabilities that a walk starting at $x$ will reach $v_1, v_2$ or $u$ first. Then, by linearity of expectation,
    \begin{equation}
    \label{eqn:proba_x_u}
        H(x) = aH(v_1) + bH(v_2) + (1-a-b)H(u) = a + (1-a-b)H(u)\implies H(u)=\frac{H(x)-a}{1-a-b}. 
    \end{equation}
    Let $A$ be the set of walks that reach $v_1$ before $u$, and let $B$ be the set of walks that never touch $u$. Then, $a=P(A)=P(A|B)P(B)=\tilde H(x) (a+b)$ since
  the set of walks in $\hcell$ that start at $x$ and never touch $u$ are in bijection with the set of walks that start at $x$ in $\tilde{C}$.  Therefore $\tilde H(x)=\frac{a}{a+b}$, and if $H(u)\leq H(x)$, 
\begin{eqnarray*}
    \frac{H(x)-a}{1-a-b}=H(u)&\leq& H(x) \implies H(x)\leq \frac{a}{a+b}\\
\implies H(x) &\leq& \tilde H(x).    
\end{eqnarray*}
\end{proof}
We can now prove the main result of the section. 
\begin{proof}[Proof of Theorem~\ref{thm:alpha_mu}]
Given a sequence of graphs, we create the sequence of harmonic functions on those graphs with boundary values on $v_1$ and $v_2$ of 1 and 0 respectively.  Let $(u_1, = v_1, u_2=w, \dots, u_{d+1} = v_2)$ be a shortest path from $v_1$ to $v_2$ in $\hcell$.  The graph sequence is built by iteratively moving along the path, and removing all the neighbours of vertices on the path, except $u_1, \dots, u_{d+1}$. In every subgraph, as $w$ is the only vertex adjacent to $v_1$, every walk starting from a different vertex must pass through $w$. Thus, the harmonic function on this graph, $H'$, satisfies $H'(v)\leq H'(w)$ for all $v\neq v_1$ in the subgraph. If $v$ is off of the path, then $v_1$, $v_2$ and $w$ will be in the same component if $v$ is deleted from the graph. Thus, in the graph with $v$ removed, the value of its harmonic function $\tilde{H}$ satisfies $H(w)\leq \tilde H(w)$ by the above lemma.

Thus, in the sequence of harmonic functions generated by the the sequence of graphs, the value on $w$ (or indeed, any point on the path) is monotonically increasing. The process terminates at a path, the smallest connected graph, since removing any vertex will separate $v_1$ and $v_2$. 

Therefore, in the general case $H(w) \leq H_L(w)=\frac{d(v_1,v_2)
  -1}{d(v_1,v_2)}$ and so $\alpha \leq d(v_1, v_2) \leq \mu$. The last
inequality is an equality if and only if $\hcell$ is a finite path graph
since $\mu$ is the total number of edges in the graph. 
\end{proof}

\section{Perspectives}
A proof of Conjecture~\ref{conj:main} with this argument framework
could be accomplished with two additional results. First, we could
prove a version of Theorem~\ref{thm:alpha_mu} for $\theta>2$. Our
arguments use the fact that $v_1$ has only one neighbour. With some
additional bookkeeping, perhaps the argument could be extended. The
core is doubly vertex transitive, and this might be useful. Next, one must show that when $X$ is not bipartite, and $G$ is
algebraic, $Q$ only has a zero at $1$. The fact that $G^N$ is rational
is a very strong characterization and tells that each singularity is
isolated. This, in combination with known information about the Julia
set may yield the desired result.  Establishing the more general result of Kr\"on and Teufl would then be reduced to showing that the only case of a constant~$\omega$ occurs with an algebraic Green's function. 

\section*{Acknowledgments} This work was supported by Canadian
National Science and Engineering Research Council (NSERC) Discovery
Grant ``Transcendence and Combinatorics'', and early readers that
provided extremely valuable feedback. 
\bibliographystyle{amsplain}
\bibliography{refs}
\end{document}